\documentclass[12pt]{amsart}
\usepackage{amssymb,amsmath}
\newtheorem*{thma}{Theorem  of Harvey-Lawson}
\newtheorem*{thmb}{Theorem of Rothstein-Sperling}

\def\4#1{\mathbb{#1}}
\def\B{\4B}

\numberwithin{equation}{section}
\def\db{\bar\partial}
\def\db*{\bar\partial^*}

\def\G{{\mathcal G}}
\def\U{{\mathcal U}}
\def\I{{\mathcal I}}
\def\cR{{\mathcal R}}
\def\W{{\mathcal W}}
\def\T{\text}

\def\1#1{\overline{#1}}
\def\2#1{\widetilde{#1}}
\def\3#1{\widehat{#1}}
\def\4#1{\mathbb{#1}}

\def\5#1{\frak{#1}}
\def\6#1{{\mathcal{#1}}}
\def\C{{\4C}}
\def\R{{\4R}}

\def\Z{{\4Z}}

\def\dib{\bar\partial}
\begin{document}
\title[Boundaries of analytic varieties]{Boundaries of analytic varieties}
\author[L.~Baracco]{Luca Baracco}
\address{Dipartimento di Matematica, Universit\`a di Padova, via 
Trieste 63, 35121 Padova, Italy}
\email{baracco@math.unipd.it}
\maketitle
\begin{abstract}
We prove that every smooth CR manifold $M\subset\subset \C^n$, of hypersurface type, has a complex strip-manifold extension in $\C^n$. If $M$ is, in addition, pseudoconvex-oriented, it is the ``exterior "boundary of the strip. In turn, the strip extends to a variety with boundary $M$ (Rothstein-Sperling Theorem); in case $M$ is contained in a pseudoconvex boundary with no complex tangencies, the variety is embedded in $\C^n$. Altogether we get: $M$ is the boundary of a variety (Harvey-Lawson Theorem); if $M$ is pseudoconvex oriented the singularities of the variety are isolated in the interior; if $M$ lies in a pseudoconvex boundary, the  variety is embedded in $\C^n$ (and is still smooth at $M$).
\end{abstract}
\def\Giialpha{\mathcal G^{i,i\alpha}}
\def\cn{{\C^n}}
\def\cnn{{\C^{n'}}}
\def\ocn{\2{\C^n}}
\def\ocnn{\2{\C^{n'}}}
\def\const{{\rm const}}
\def\rk{{\rm rank\,}}
\def\id{{\sf id}}
\def\aut{{\sf aut}}
\def\Aut{{\sf Aut}}
\def\CR{{\rm CR}}
\def\GL{{\sf GL}}
\def\Re{{\sf Re}\,}
\def\Im{{\sf Im}\,}
\def\codim{{\rm codim}}
\def\crd{\dim_{{\rm CR}}}
\def\crc{{\rm codim_{CR}}}
\def\phi{\varphi}
\def\eps{\varepsilon}
\def\d{\partial}
\def\a{\alpha}
\def\b{\beta}
\def\g{\gamma}
\def\D{\Delta}
\def\Om{\Omega}
\def\k{\kappa}
\def\l{\lambda}
\def\L{\Lambda}
\def\z{{\bar z}}
\def\w{{\bar w}}
\def\Z{{\1Z}}
\def\t{{\tau}}
\def\th{\theta}
\emergencystretch15pt
\frenchspacing
\newtheorem{Thm}{Theorem}[section]
\newtheorem{Cor}[Thm]{Corollary}
\newtheorem{Pro}[Thm]{Proposition}
\newtheorem{Lem}[Thm]{Lemma}
\theoremstyle{definition}\newtheorem{Def}[Thm]{Definition}
\theoremstyle{remark}
\newtheorem{Rem}[Thm]{Remark}
\newtheorem{Exa}[Thm]{Example}
\newtheorem{Exs}[Thm]{Examples}
\def\Label#1{\label{#1}}
\def\bl{\begin{Lem}}
\def\el{\end{Lem}}
\def\bp{\begin{Pro}}
\def\ep{\end{Pro}}
\def\bt{\begin{Thm}}
\def\et{\end{Thm}}
\def\bc{\begin{Cor}}
\def\ec{\end{Cor}}
\def\bd{\begin{Def}}
\def\ed{\end{Def}}
\def\br{\begin{Rem}}
\def\er{\end{Rem}}
\def\be{\begin{Exa}}
\def\ee{\end{Exa}}
\def\bpf{\begin{proof}}
\def\epf{\end{proof}}
\def\ben{\begin{enumerate}}
\def\een{\end{enumerate}}
\def\dotgamma{\Gamma}
\def\dothatgamma{ {\hat\Gamma}}

\def\simto{\overset\sim\to\to}
\def\1alpha{[\frac1\alpha]}
\def\T{\text}
\def\R{{\Bbb R}}
\def\C{{\Bbb C}}
\def\Z{{\Bbb Z}}
\def\Fialpha{{\mathcal F^{i,\alpha}}}
\def\Fiialpha{{\mathcal F^{i,i\alpha}}}
\def\Figamma{{\mathcal F^{i,\gamma}}}
\def\Real{\Re}
%
%
%
\section{Boundaries of analytic manifolds and analytic varieties }
\Label{s1}

The purpose of this paper is to give a comprehensive description of the  boundaries of the analytic varieties from the point of view of the CR Geometry.
In particular, to offer a new, simple proof of the Harvey-Lawson's Theorem and to improve its conclusion when the boundary is weakly pseudoconvex. 
 It starts by a construction of families of analytic discs which yields
the existence of a strip-manifold extension of every smooth CR submanifold, without boundary, of hypersurface type  $M$ in $\C^n$, compact, connected and oriented (Theorem~\ref{t1.1} point (i)); moreover,
when $M$ is pseudoconvex, it is the exterior, pseudoconvex component of the boundary of the strip (point (ii)).
Next, there is discussed the extension from the strip to its ``interior", first by recalling
 the Rothstein-Sperling Theorem \cite{RS65};
 for the selfcontainedness of the paper, there is presented in Section~\ref{s3} a sketch of the proof; it follows the proof by Yau in \cite{Y81} which is in turn inspired to the one by Siu in \cite{S74} but  also uses in an essential way Proposition~\ref{p2.1} below.
Next, it is proved that  if the strip is contained in a pseudoconvex domain to whose boundary it is transversal, then  the extension is fully embedded in $\C^n$ (Theorem~\ref{t1.2}).
Combination of Theorem~\ref{t1.1} with the Rothstein-Sperling Theorem yields the extension from $M$ to a variety, possibly singular at $M$; this is  the celebrated Harvey-Lawson Theorem of \cite{HL75}. When one starts from $M$ pseudoconvex, the  extension variety $\W$  encounters  the non singular strip of $\C^n$  before approaching  its boundary $M$. Thus $\W$ is smooth in a neighborhood of $M$ and one gets back the result by the author in \cite{B12} (Theorem~\ref{t1.3} which extends Theorem~10.4 of \cite{HL75}). 
 This, in combination with the invariant estimates  of \cite{K86}, were used in \cite{B12} to solve  a conjecture by Kohn: if $M\subset\subset \C^n$ is pseudoconvex, then $\dib_b$ has closed range.
 Finally, combination of Theorem~\ref{t1.1} and \ref{t1.2} yields that if $M$ is contained in a pseudoconvex boundary, and it is not complex tangential, then the variety that it bounds is embedded in $\C^n$ and still has only isolated singularities in the interior (Theorem~\ref{t1.4}).
 Note that Luk and Yau have given in \cite{LY98} a simple example of a strongly pseudoconvex manifold $M$ which bounds a variety of $\C^n$ which is singular at $M$ in apparent contrast to \cite{HL75} Theorem~10.4 and to our Theorem~\ref{t1.3}. According to the Fefferman's Editor Note \cite{F} this is just a   conflict of terminology: the variety $\W$ with boundary $M$ is singular as ``embedded" but smooth as ``immersed". That means that $\W$ is equipped with a smooth mapping $\pi:\,\W\to \C^n$, holomorphic in $\W\setminus M$ which is an embedding only in a neighborhood of $M$ but not  in the whole $\W$. This explains Theorem~\ref{t1.3}. Next, what we wish to point out in Theorem~\ref{t1.4} is that, when one starts from $M$ contained in a pseudoconvex boundary with no complex tangencies, there is no possible   conflict of terminology at all: the variety sits in $\C^n$ and it is smooth at $M$.

Here is our first result, whose second item (ii) is already contained in \cite{B12}; it is a consequence of the theory of minimality in the CR geometry.

\bt
\Label{t1.1}
 \hskip0.2cm (i) Let $M\subset\subset\C^n$ be a smooth, compact, CR manifold of hypersurface type. Then there is a complex strip-manifold $Y$ whose closure contains $M$ and which is smooth up to $M$. 

\noindent
(ii) When $M$ is, in addition, pseudoconvex-oriented, then $M$ is the exterior, that is pseudoconvex, boundary of $Y$.
\et

Before starting the proof we need a few preliminaries about CR Geometry. 
Let $M$ be a CR manifold of hypersurface type in $\C^n$, denote by $TM$ and $T^\C M:=TM\cap JTM$ the tangent and complex tangent bundle respectively, and let $\pm T$, respectively $\pm \omega$, be  purely imaginary generators of the quotient bundle $\frac{TM}{T^\C M}$, respectively of the dual bundle $T^\C M^\perp$. We say that $M$ is pseudoconvex when, for a choice of the sign, $+d\omega|_{T^\C M}\ge0$.
A CR curve $\gamma$ on $M$ is a real curve such that $T\gamma\subset T^\C M$ and a CR orbit is the union of all piecewise smooth CR curves issued from a point of $M$. According to Sussmann's theorem (cf. \cite{MP06}) the orbit has the structure of an immersed variety of $\C^n$. 
\vskip0.3cm
\noindent
{\it Proof of Theorem~\ref{t1.1}. }\hskip0.2cm(Cf. \cite{B12}). Our proof begins by showing that, if $M\subset\subset \C^n$ satisfies the hypotheses of the Theorem, then it consists of a single CR orbit.
We   point out that this result is classical for a hypersurface, the boundary of a domain of $\C^n$ (cf. \cite{J95} and \cite{MP06} Lemma 4.18). First,
since $M$ is compact, the inductive  family of the CR-invariant closed non-empty subsets has a non-empty minimal $S$; this is the closure of the orbit  $\mathcal O$ through whatever point of $S$. We claim that $\mathcal O$ is open. Otherwise, $\mathcal O$ and indeed $S$ itself by density, cannot contain any point of minimality (in the sense of \cite{T88}).  Since $M$ is  hypersurface type, without boundary, $S$ is  foliated by complex 1-codimensinal leaves, without boundary, each one being dense in $S$.  Now, each coordinate function $z_i,\,\,i=1,...,n$ is constant in each leaf where $|z_i||_S$ achieves a maximum, hence in the whole $S$ by density.
This yields contradiction; thus $\mathcal O$ is open and $S\setminus \mathcal O$  is empty since, otherwise, $\bar{\mathcal O}$ would contradict  the minimality of $S$. It follows that $\mathcal O$ is both open and closed and therefore it coincides with $M$ by connectedness.

We are ready to carry out the proof.
First, take a point $p_o$ of local minimality, that is, a point through which there passes no complex submanifold $S\subset M$; this set is certainly non-empty. 
Take  a local patch $M_o$ at $p_o$  in which the projection $\pi_{p_o}:\C^n\to T_{p_o}M+iT_{p_o}M$ induces a diffeomorphism between $M_o$ and $\pi_{p_o}(M_o)$. Since $\pi_{p_o}(M_o)$ is a  piece of a minimal  hypersurface, then $(\pi_{p_o}|_{M_o})^{-1}$ extends holomorphically to either side $\pi_{p_o}(M_o)^\pm$ by  \cite{T88}, and parametrizes a one-sided complex manifold $Y_o$ which has a neighborhood of $p_o$ in $M_o$ as its boundary. 
We rephrase this conclusion by saying that $M$ extends to either  direction $\pm J T$ at any point of local minimality.
If $p$ is now a general point of $M$, connect $p$ to a point $p_o$ of local minimality by a piecewise smooth CR curve $\gamma$.  By uniqueness of holomorphic functions having the same trace on a real hypersurface, one-sided complex  neighborhoods glue together into a complex neighborhood of a maximal open  arc $\gamma_o\subset\gamma$, which starts  from $p_o$ and ends at  $p_1$. If $p_1\neq p$, then $\pi_{p_1}(M)$ is not minimal and therefore there exists a complex hypersurface $S\subset \pi_{p_1}(M)$ which contains $\pi_{p_1}(p_1)$. Also, $\pi_{p_1}(\gamma)$ being a CR curve, it must belong to $S$. Thus,
 extension of $(\pi_{p_1}|_M)^{-1}$ to $\pm\pi'_{p_1}(JT)$ propagates from $\gamma_o$ along $S$ beyond $\pi_{p_1}(p_1)$ by Hanges-Treves Theorem \cite{HT83}, a contradiction.

\hskip11cm $\Box$

\noindent
 We wish to continue the strip to a variety with boundary.
 \begin{thmb}
 \Label{thmb}
 There is a complex variety $\mathcal W$ which extends the strip $Y$.
 \end{thmb}
 In the terminology of \cite{RS65}, $\W$ is a variety whose boundary is a cycle of $Y$ homotopic to $M$.
A sketch of the proof follows in Section~\ref{s3}.
In a special case we have a better conclusion.
\bt
\Label{t1.2}
Assume that $Y$ is contained in a pseudoconvex domain $\Omega\subset\C^n$ and is transversal to its boundary. Then there is a complex variety $W\subset\C^n$ whose boundary is the pseudoconvex boundary of the strip and which has only isolated singularities in the interior.
\et
The point here is that $W$, differently from $\W$, is embedded in $\C^n$. The proof follows in Section~\ref{s2}.
\begin{thma}
\Label{thma}
 Let $M\subset\subset \C^n$ be a compact, smooth, oriented, hypersurface type submanifold. Then $M$ is the boundary of a complex variety $\W$. 
 \end{thma}
 \bpf
 This is a combination of Theorem~\ref{t1.1} (i) and the Rothstein-Sperling Theorem.
 
 \epf
 In general, $\mathcal W$ is not smooth in a neighborhood of $M$. According to Harvey and Lawson, this is true when $M$ is strongly pseudoconvex; we show here that weak pseudoconvexity is indeed sufficient. 
 \bt
 \Label{t1.3}
 (cf. \cite {B12})
Suppose that, in addition to the hypotheses of Theorem~\ref{thma}, the manifold $M\subset\subset\C^n$  is also pseudoconvex-oriented. Then $M$ is the boundary of a complex variety $\W$ smooth in a neighborhood of $M$.
\et
\bpf
Since $M$ is pseudoconvex-oriented, then by Theorem~\ref{t1.1} (ii), we know that $M$ is the exterior, pseudoconvex boundary of the smooth strip; therefore the singularities of $\mathcal W$ cannot cluster at $M$.

\epf
\bt
\Label{t1.4}
Assume, in addition to the hypotheses of the Harvey-Lawson Theorem, that $M$ is contained in the boundary of a pseudoconvex domain $\Omega\subset\C^n$ to which it is not complex tangential. Then $M$ is the boundary of a complex  variety $W\subset\C^n$ with  only isolated singularities in the interior.
\et
\bpf
First note that any disc attached to $M$ is forced to lie in $\overline{\Omega}$ by pseudoconvexity. It follows that the whole strip $Y$ constructed in Theorem~\ref{t1.1}  is contained in $\overline{\Omega}$. We claim that it  is transversal to $\partial\Omega$ at every point of $M$. In fact if $Y$ were tangential to $\partial\Omega$ at a point $p_o\in M$, then all discs of $Y$ passing through $p_o$  would be contained in $\partial\Omega$. Therefore an open subset $U$ of $Y$ would be contained in $\partial\Omega$. 
By taking discs through different points of $U$, we see that the subset $U\subset Y$ of points belonging to $\partial\Omega$ can be enlarged, step by step, to the full $Y$ which implies $TM\subset T^\C \partial\Omega$, a contradiction. Thus,  the strip generated by $M$ is transversal to $\partial\Omega$ at every point, and $M$ is its pseudoconvex-oriented boundary. We then apply Theorem~\ref{t1.2} to extend $Y $ to a variety $W\subset \C^n$ with smooth boundary $M$.

\epf

\section{Analytic extension from the boundary of a pseudoconvex domain - Proof of Theorem~\ref{t1.2}.} 
\Label{s2}
In this Section, we show  how to continue the strip produced by Theorem~\ref{t1.1} to a variety with boundary
 $M$ when $M$ sits in a pseudoconvex hypersurface.
 For this,
we recall a few facts about analytic sets. Let $D$ be a domain of $\C^n$ and $A$ a subset of $D$.
\bd The set $A$ is said to be analytic if for all points $z\in D$, there exist a neighborhood $V$ of $z$ and a system of holomorphic functions $h_1,...,h_d \in \mathcal{O} (V)$ such that $V\cap A= \{ w\in V| h_1(z)=0,...=h_d(z)=0 \}$. The set $A$ is said to be locally analytic if the  above property  only holds for the points $z$ of $A$ instead of  the full $D$.
\ed
Analytic sets are locally the common zero-sets of holomorphic functions, in particular, they  are locally closed. Complex submanifolds are locally analytic sets.
\bp[Extension through a pseudoconcave boundary] 
\Label{p2.1}
\label{ext}
Let $D\subset\C^n$ be a strictly pseudoconcave domain and $A$ an analytic set of $D$ in a neighborhood of a point $p\in \partial D$. Then there is an analytic extension of $A$ from $D$ across $\partial D$.
\ep
\bpf
It is not restrictive to assume  $p\in \bar A$. By a coordinate change we may further suppose that $p=0$, that $ T^\C_p\partial D$ is defined by $z_n=0$, and that $D$ is strictly concave. We denote by $A'$ the intersection of $A$ with the plane $z_n=0$. By Stein-Remmert's Theorem, $A'$ extends across $0$. Decompose the variable as $z=(z',z'',z_n)$ for $z'=(z_1,...,z_{n-k})$ and  $z''=(z_{n-k+1},...,z_{n-1})$, denote by $\pi_1:\,\C^n_z\to \C^{n-k}_{z'}$ the projection, and by $\B^{n-k}_\delta$, $\B^{k-1}_\delta$ and $\B^{k}_\delta$ the balls of radius $\delta$ in $\C^{n-k}_{z'}, \C^{k-1}_{z''}$ and $\C^{k}_{(z'',z_n)}$ respectively.    
Suppose that $A'$ is normalized in the $z'$ directions; hence, for suitable $\delta_1$, the projection $\pi_1:\,A'\to \B^{n-k-1}_{\delta_1}$ has discrete fibers. (To see it, we just have to think of $A'$ as the zero set of  a system of Weierstrass polynomials in  $z'$.) Thus we can find $\delta_2<\delta_3$ such that 
\begin{equation}
\Label{nova}
\Big ((\B^{n-k}_{\delta_3}\setminus \B^{n-k}_{\delta_2})\times (\B^{k-1}_{\delta_1}\times \{0\})\Big)\cap A=\emptyset.
\end{equation}
By continuity, \eqref{nova} remains true when we replace $\B^{k-1}_{\delta_1}\times \{0\}$ by $\B^{k}_{\delta_1}$. Let us define $V_o:=\{(z'',z_n):\,\B^{n-k}_{\delta_3}\cap \pi_1^{-1}(z'',z_n)\subset D\}$; this is a strictly concave set. Recall that $A\subset D$; by \eqref{nova}, applied for $(z'',z_n)\in V_o$, we have that $\pi_1:\,A\to V_o$ has compact, hence finite fibers in $\B^{n-k}_{\delta_3}$. It follows that $A\cap \Big(\B^{n-k}_{\delta_3}\times V_o\Big)$ is globally defined over $\B^{n-k}_{\delta_3}\times V_o$ as the zero set of $n-k$ Weierstrass polynomials in $z'$ with analytic coefficients in $(z'',z_n)$, that is 
$$
P^1_{(z'',z_n)}(z_1)=0,...,P^{n-k}_{(z'',z_n)}(z_{n-k})=0.
$$
The analytic continuation of the coefficients from $V_o$ to its convex hull implies the extension of their zero set $A$ from $D$ to $\B^{n-k}_{\delta_3}\times V_o$.

\epf

\br In general, if the hypersurface has less than $n-1$ but more than  $n-k+1$ negative eigenvalues, then $A$ can be
extended provided that $\dim(A)\ge k$. In fact, under this assumption, the set $V_o$, which is  the intersection of $D$ with a $n-k$ complex plane, has at least one negative Levi eigenvalue and this suffices for analytic extension of holomorphic functions which yields the conclusion.                      
\er
\vskip0.3cm \noindent {\it Proof of Theorem~\ref{t1.2}.} \hskip0.2cm
We have to extend the strip $Y$  as an analytic set from a neighborhood of $\partial\Omega$ to the full $\Omega$; we do it through a family of  pseudoconcave hypersurfaces which shrink to a point.
Since $Y$ is transversal to $\partial\Omega$, then, for a suitable $\delta$, it is an anlytic set in $\Omega\setminus\Omega_\delta$ for $\Omega_\delta=\{z\in\Omega:\T{dist}(z,\partial\Omega)>\delta\}$.
Assume $0\in\Omega$, and let $S_t$, $t\in\R^+$, be the hypersurface $S_t:=\partial(\B_t\cap \Omega)$.
We have
\begin{itemize}
\item[(i)] Points of $S_t$ in $\partial\Omega$, are points  where $Y$ is already extended to the interior of $\B_t\cap \Omega$ in a neighborhood of $S_t$.
\item[(ii)] At points of $S_t$ in $\partial\B_t$, extension is a consequence of Proposition~\ref{p2.1}.
\end{itemize}
Letting $t\searrow0$, we see that $Y$ extends as an analytic set in a domain sheeted over the whole of  $\Omega\setminus\{0\}$. But, $\Omega$ being pseudoconvex, this domain is in fact $\Omega\setminus\{0\}$ and eventually $\Omega$ itself (by
trivial extension to the isolated point  $0$). 
 \hskip11cm

\section{Appendix - Sketch of proof of the Rothstein-Sperling Theorem.}
\Label{s3}

We give  an outline of the proof of the Rothstein-Sperling Theorem. The most part is inspired to  \cite{Y81} but also uses  Proposition~\ref{p2.1} above. The idea is similar to the one we used in Theorem \ref{t1.2}, namely the extension of an analytic set through a family of pseudoconcave hypersurfaces. We are no longer assuming that the strip $Y$ is transversal to a pseudoconvex boundary and thus we can no longer use the family of hypersurfaces which shrink from the annulus to the origin. Instead, we reach the conclusion by extending through a family of spheric shells. In this case, since the intersection of a sphere with the analytic set happens to have several components, there are many extension surface pieces, one for each component. However, these will eventually  merge to form the desired variety.  But this procedure will not take place in $\C^n$, where monodromy may fail but, instead, in a space $\G$ over $\C^n$.
\bd
The space $\G$ of irreducible analytic germs is defined by
$$
\mathcal{G}=\{ (p, X_p):\,\T{$p$ is a point of $ \C^n$ and $X_p$ a germ of an irreducible analytic set at } p\}.
$$
\ed
 We put a topology on $\G$. The system of open neighborhoods of a point $(p,X_p)$ is obtained from the system of the neighborhoods $U_p\subset\C^n$  of $p$  on which $X$ is globally defined and irreducible by putting $U_{(p,X_p)}=\{(q,X_q):\,q\in U_p\}$.

\bd A subset $A$ of $\G$ is said to be a local analytic set (or a surface piece) over $\C^n$  if for any point $(p,X_p)\in A$ there exists a neighborhood $U_{(p,X_p)}$ such that $\pi$ induces a homeomorphism between $A\cap U_{(p,X_p)}$ and an analytic set in $\pi(U_{(p,X_p)})$.  
\ed
The strip $Y$ being regular, it takes a trivial lift to $\G$. For  $\phi_r=-|z|^2+r^2$, if  $R_r, I_r$ and $U_r$ are the sets defined by $\phi_r=0$, $\phi_r>0$ and $\phi_r<0$ respectively, we denote by $\cR_r$, $\I_r$ and $\U_r$ their respective preimages under $\pi$.
\bd A cycle $Z$ in $Y$ is a piecewise $C^\omega$-hypersurface in $Y$ with the property that there exists a relatively compact open neighborhood $Y_1$  of $Z$ in $Y$ such that $Y_1\setminus Z$ has two connected components $Y^+_1$ and $Y^-_1$ with $Y_1=Y_1^+\cup Y_1^-\cup Z$ and $Z= \overline{Y_1^+} \cap \overline{Y_1^-}$.\ed  
Remark that
the strip $Y$ of Theorem~\ref{t1.1}  contains a 1-parameter family of real analytic cycles homotopic to the initial manifold $M$.
\vskip0.3cm
\noindent
{\it Sketch of proof of the Rothstein-Sperling Theorem.}\hskip0.2cm
What we prove is that there is a complex variety of $\G$ with boundary $Z$ which extends $Y_1^+$. 
We begin by noticing that,
 $Z$ and $\cR_r$  being real analytic, then $Z\cap\cR_r$  has a finite number of components $\{k\}$ that we call $K_r$-arcs. Likewise,  $Z\cap \U_r$ decomposes into a finite family of components $\{\alpha\}$ that we call  $A_r$-arcs.
The $A_r$-arcs are relatively open while the $K_r$-arcs are closed. The boundary of an $A_r$-arc contains points of $K_r$-arcs, but not all $K_r$-arcs have points belonging to the boundary of an $A_r$-arc. We call $T_r$-arcs  and denote by $\{\beta\}$ this particular  type of $K_r$-arcs. It is readily seen that $T_r$-arcs correspond to points where the spheric surface begins to touch $Z$ as $r$ decreases. Let $r_o=\inf\{r>0 | Z\cap \U_r =\emptyset\}$; then $Z\cap \cR_{r_o}$ consists only of $T_{r_o}$-arcs  $\{\beta\} $. For every $\beta$, choose a small neighborhood $U=U_\beta$ in $Y_1$. This is divided in two sides by $Z$ that we denote by $U^+$ and $U^-$; again, these depend on $\beta$. Choose the side whose interior lies entirely in $\I_{r_o}$ and refer to it as the positive side. Then for $r$ slightly smaller than $r_o$ we have that every $A_{r}$-arc lies in the interior of some $U$. For these $r$, we have that for every $A_r$-arc $\alpha$, there exists  a piece of surface $f$ whose boundary contains $\alpha$ and which satisfies $\partial f\setminus \alpha \subset \cR_r$ and $ f\subset \U_r$.  
Thus, for $r$ close to $r_o$, the family $\{\alpha\}$ satisfies what we call the ``property $(E_r)$". This consists in the existence of a family of pieces of surfaces $\{f\}$ such that $\alpha \subset \partial f$, $\partial f\setminus\alpha\subset\cR_r$, and $\alpha\not\subset\partial f'$ (if $f'$ is the piece of surface which satisfies $f'\supset\alpha'$ for a different $\alpha'$).  Moreover, each $f$ is requested to lie in $\U_r$  locally on one side of $Z$ in a neighborhood of $\alpha$.

 We define an $\hat A_r$-arc as a connected component of $Z\cap(\U_r\cup\cR_r)$. Every $\hat A_r$-arc is the union of 
 some $A_r$-and $K_r$-arc; assume, for instance,   $\hat{\alpha}=\alpha_1\cup\alpha_2\cup k$. Suppose that $(E_r)$ holds; then there are  
 $f_1$ and $f_2$ such that $\partial f_1=\alpha_1$ and $\partial f_2=\alpha_2$. One can readily check that  $f_1$ and $f_2$ lie in the same 
 side of $Z$ (cf. \cite{Y81} Lemma 6.22). Let  $V$ be a small neighborhood of $k$ contained in $U$. Since, by Proposition~
 \ref{p2.1}, $f_1$ and $f_2$ extend  through $\cR_r$ to two  pieces of surfaces $\tilde f_1$ and $\tilde f_2$, then by taking the union of $
 \tilde f_1$, $\tilde f_2$ and $V$, we get a piece of surface $\tilde f$.
 What we proved for $\hat\alpha=\alpha_1\cup\alpha_2\cup k$ is true for any $\hat A_r$-arc. 
 With this preparation, we can prove $(E_{r'})$ for $r'$ close to $r$ with $r'<r$. In fact, 
  for those $\alpha'$ having non-empty intersection with $\U_r$, i.e. $\alpha\subset \alpha'$ for some $A_r$-arc $\alpha$, we have that the piece of surface $\tilde f \cap \U_{r'}$ is bounded by $\alpha' $. The other $A_{r'}$-arcs arise from $T_r$-arcs and for these, if $r'$ is sufficiently close to $r$, there exist pieces of surfaces $f $ bounded by them (cf. the remark before the definition of property $(E_r)$). Altogether, if $(E_r)$ holds for $r$, it also holds  for some $r'$ smaller than $r$. Since it trivially holds for $r$ next to initial value $r_o$, then it holds for any $r>0$. 
  
 The different extensions obtained for different $ A_r$-arcs and different $r$ glue up to a single $\tilde f$ whose boundary has,  when $r$ is small, two componets, the cycle $Z$ itself which represents the unique $A_r$-arc, and another component on $\partial \B_r$. Extending this latter across $\partial \B_r$ by means of Proposition~\ref{p2.1}, we extend $\tilde f$  to $\G\setminus\pi^{-1}(0)$. Changing the center of the system of spheres $\B_r$, we extend $\tilde f$ to the whole $\G$.
  
\hskip11cm $\Box$
\vskip0.3cm
\noindent
{\it Acnowledgments.}\hskip0.2cm  I am greatly indebted to Alexander E. Tumanov for several useful discussions
 and, especially, for directing my attention to the crucial role played
 in Complex Geometry by CR minimality.


\begin{thebibliography}{KN06}
\bibitem{B12}{\bf L. Baracco}---The range of the tangential Cauchy-Riemann system to a CR embedded manifold, {\em Invent. Math.} (2012)

\bibitem{F}{\bf C. Fefferman}---Editor's note on papers by F. R. Harvey and H. B. Lawson, Jr.: "On boundaries of complex analytic varieties. I'' [Ann. of Math. (2) 102 (1975), no. 2, 223--290; MR0425173 (54 13130)] and by H. S. Luk and S. S.-T. Yau: "Counterexample to boundary regularity of a strongly pseudoconvex CR submanifold: an addendum to the paper of Harvey-Lawson'' [ibid. 148 (1998), no. 3, 1153--1154; MR1670081 (99j:32006)],
 {\em Ann. of Math.} {\bf 151}  (2)  (2000),  875

\bibitem{FG01}{\bf K. Fritzche and H. Grauert}---From holomorphic functions to complex manifolds,{\em Springer Graduate texts in Mathematics} {\bf 213} ,2001
\bibitem{HT83}{\bf J. Hanges and F. Treves}---Propagation of holomorphic extandibility of CR functions, {\em Math. Ann.} {\bf 263} n. 2 (1983), 157--177
\bibitem{HL75} {\bf F.L. Harvey and H.B. Lawson}---On boundaries of analytic varieties, I, {\em Ann. of Math.} {\bf 102} (1975), 223--290
\bibitem{J95}{ \bf B. Joricke}---Some remarks concerning holomorphically convex hulls and envelopes of holomorphy, {\em Math. Z.} {\bf 218} n. 1 (1995), 143--157 
 
\bibitem{K86} {\bf J. J. Kohn}---The range of the tangential Cauchy-Riemann operator, {\it Duke Math.
J\/}.\ {\bf 53} (1986), 525--545

\bibitem{LY98} {\bf H. S. Luk and S. T. Yau}---Counterexample to boundary regularity of a strongly pseudoconvex CR submanifold: an addendum to the paper of Harvey-Lawson, {\it Ann. of Math. (2)} {\bf 148} (1998), no. 3, 1153???1154 

\bibitem{MP06}{\bf J. Merker and E. Porten}---Holomorphic extension of CR functions, envelopes of holomorphy, and removable singularities, {\em IMRS Inst. Ter. Survey} (2006)
\bibitem{RS65}{\bf W. Rothstein and H. Sperling}---Einsetzen analytischer Fl\"achenst\"ucke in Zyklen auf komplexen R\"aumen. (German) 1966 Festschr. {\it Ged\"achtnisfeier K. Weierstrass} 531???-554 Westdeutscher Verlag, Cologne

\bibitem{S74} {\bf Y.-T. Siu}---Techiques of extension of analytic objects, {\it Lect. Notes in Pure and Appl. Math.} , {\bf 8} Marcel Dekker Inc., New York, (1974)

\bibitem{T88} {\bf A. Tumanov}---Extending CR functions on a manifold of finite type over a wedge, {\em Mat. Sb.} {\bf 136} (1988), 129--140

\bibitem{T94} {\bf A. Tumanov}---Connection and propagation of
  analyticity of CR functions, {\em Duke Math. J.} {\bf 71}, n.\ 1 (1994), 1--24
\bibitem{Y81}{\bf S. Yau}---Kohn-Rossi cohomology and its application to the complex Plateau problem I, {\em Ann. of Math.} {\bf 113} (1981), 67--110



 

 

















 

 













 
 
 
 

 


 
\end{thebibliography}
\end{document}